# DISTRIBUTION OF POWERS MODULO 1 AND RELATED TOPICS

MIGUEL A. LERMA

Abstract. This is a review of several results related to distribution of powers and combinations of powers modulo 1. We include a proof that given any sequence of real numbers $\theta_n$, it is possible to get an $\alpha$ (given $\lambda \neq 0$), or a $\lambda$ (given $\alpha > 1$) such that $\lambda\,\alpha^n$ is close to $\theta_n$ modulo 1. We also prove that in a number field, if a combination of powers $\lambda_1\,\alpha_1^n + \cdots + \lambda_m\,\alpha_m^n$ has bounded $v$-adic absolute value (where $v$ is any non-Archimedean place) for $n \geq n_0$, then the $\alpha_i$'s are $v$-adic algebraic integers. Finally we present several open problems and topics for further research.

## 1. Introduction

The study of the behaviour of sequences of the form $\alpha^n$ modulo 1 has some interesting connections with subjects such as Waring's problem. Let $g$ be the function:

$$g(k) = \min\{\, s \in \mathbb{N} \,:\, a = n_1^k + \cdots + n_s^k \text{ for all } a \in \mathbb{N}\,\} \tag{1.1}$$

and $\|x\|$ = distance from x to the nearest integer. Then it is well known that for $k \geq 5$, if

$$\|(3/2)^k\| \;>\; (3/4)^k \tag{1.2}$$

then

$$g(k) \;=\; 2^k + [(3/2)^k] - 2 \tag{1.3}$$

so that the rate to which $(3/2)^k \pmod 1$ accumulates near 0 determines the value of $g(k)$.

In a slightly more general setting, the study of sequences of the form $\lambda\,\alpha^n$ also has interesting applications. As an example, it is known ([5], theorem 8.1) that a number $\lambda$ is normal to the base $b$, i.e., all its finite sequences of k digits to base b occur with the same relative frequency $1/b^k$, if and only if the sequence $\lambda\,b^n$ is u.d. mod 1 (uniformly distributed modulo 1, see definition 3.2 below), so that techniques from the theory of distribution modulo 1 apply to this problem.

*Date.* June 19, 1995.

I am grateful to Prof. Jeffrey Vaaler for his support and advice. I am also grateful to Prof. David Boyd and Prof. Alf van der Poorten for their valuable help.

The study of slowly growing sequences such as $\{n\omega\}_{n=1}^{\infty}$ is not hard, but when the sequence grows very rapidly, as with $x_n = \alpha^n$ for $\alpha > 1$, the fractional part of $x_n$ becomes almost negligible compared to its integer part. So, special techniques have had to be developed to tackle the problem in that case.

Several results related to distribution of powers modulo 1 are presented below.

## 2. Overview of results

A natural question is if a sequence of the form $\alpha^n$, or more generally, of the form $\lambda\,\alpha^n$ for some fixed $\lambda \neq 0$, is u.d. mod 1. Koksma's metric theorem 3.4 shows that in fact, $\lambda\,\alpha^n$ is u.d. mod 1 for almost every $\alpha > 1$, i.e., the exceptional set of $\alpha$'s for which $\lambda\,\alpha^n$ is not u.d. mod 1 has Lebesge measure zero. It is surprising, however, that there is no known concrete example of a real number $\alpha > 1$ for which $\alpha^n$ is u.d. mod 1; only members of the exceptional set are known. For instance, the following classes of numbers are known to be in the exceptional set:

(1) Integers ($> 1$), since $\alpha^n = 0 \pmod 1$ for every $\alpha \in \mathbb{Z}$.

(2) Pisot-Vijayaraghavan (or P.V., or Thue) numbers. A P.V. number is a real algebraic integer $\alpha > 1$ whose conjugates lie inside the open unit disc $\{\, z \in \mathbb{C} : |z| < 1 \,\}$ (the rational integers greater than 1 are P.V. numbers). If $\alpha$ is a P.V number then $\lim_{n\to\infty} \alpha^n = 0 \pmod 1$ geometrically ([12], p. 3; [1], theorem 5.3.1).

(3) Salem numbers. A Salem number is a real algebraic integer $\alpha > 1$ whose conjugates lie all in the closed unit disc $\{\, z \in \mathbb{C} : |z| \leq 1 \,\}$, and at least one of them is in the border of the disc (actually it can be readily seen that all of them except one will be in the border). If $\alpha$ is a Salem number, then $\{\alpha^n\}_{n=1}^{\infty}$ is dense modulo 1, i.e., the fractional parts of $\alpha^n$ are dense in the interval $[0, 1)$, but it is not u.d. mod 1 ([12], p. 33; [1], theorem 5.3.2).

All those examples are algebraic numbers, so it is natural to ask if there is any transcendental number $\alpha$ such that $\alpha^n$ is not u.d. mod 1. A result in this direction is the following theorem of Boyd ([2]):

**Theorem 2.1.** *Let $A$, $B$ be real numbers with $3 < A < B$, and let $a_0$ be an integer satisfying $a_0 > (A+1)\,(A-1)^{-1}\,(B-A)^{-1}$. Then there is an uncountable set $S \subset [A, B]$, such that for each $\alpha \in S$, there is a real number $\lambda = \lambda(\alpha) > 0$ for which*

$$\|\lambda\,\alpha^n\| \;\leq\; (A-1)^{-1}\,(\alpha-1)^{-1} \qquad \text{for } n = 0, 1, \dots \tag{2.1}$$

*The integer $a_0$ is the nearest integer to $\lambda(\alpha)$ for all $\alpha \in S$.*

Since $S$ is uncountable, it will contain transcendental numbers. On the other hand, (2.1) shows that $\lambda\,\alpha^n$ will be in a small interval around zero modulo 1, so it cannot be u.d. mod 1.

The idea behind the proof of theorem 2.1 is to get a sequence of positive integers $a_n$ that will play the role of nearest integers to $\lambda\,\alpha^n$. Then $\lambda$ and $\alpha$ will be obtained as the following limits:

$$\alpha \;=\; \lim_{n\to\infty} \frac{a_{n+1}}{a_n} \qquad \text{and} \qquad \lambda \;=\; \lim_{n\to\infty} a_n \alpha^{-n} \tag{2.2}$$

Here, $a_0$ is given in the hypothesis of the theorem, $a_1$ will be any integer such that

$$a_0\,A + (A-1)^{-1} < a_1 < a_0\,B - (A-1)^{-1} \tag{2.3}$$

and for $n \geq 1$:

$$a_{n+1} \;=\; [a_n^2/a_{n-1}] + f(n) \tag{2.4}$$

where $[x]$ is the integer part of x, and $f \in \mathfrak{J}$ = the set of functions $\mathbb{Z}^+ \to \{0,1\}$. Since $\mathfrak{J}$ is uncountable, and each function $f \in \mathfrak{J}$ gives a different $\alpha$, the set of $\alpha$'s that can be found this way is uncountable. Furthermore, it can be proved that

$$|a_n - \lambda\,\alpha^n| \;\leq\; (A-1)^{-1}\,(\alpha-1)^{-1} \qquad \text{for } n = 0, 1, \ldots \tag{2.5}$$

which gives (2.1) and ensures that $a_n$ is, in fact, the nearest integer to $\lambda\,\alpha^n$.

However, this result does not say anything about $\lambda = 1$ or any fixed value of $\lambda \neq 0$.

The answer to the question is actually positive, as shown in theorem 3.5. Given $\lambda > 0$, and given any sequence $\theta_n$, we can get an $\alpha$ such that $\lambda\,\alpha^n$ is close to $\theta_n$ mod 1. The idea is to start with an "initial" value $\alpha_1 > A > 2$ (A large) such that $\lambda\,\alpha_1 = \theta_1 \pmod 1$. Next, slightly perturb the value of $\alpha_1$, i.e., find a slightly greater $\alpha_2 \geq \alpha_1$, such that now $\lambda\,\alpha^2 = \theta_2$. After doing so, we will not have $\lambda\,\alpha_2 = \theta_1$ (mod 1) any more, but the difference $\|\lambda\,\alpha_2 - \theta_1\|$ can be made less than $1/\alpha_1 < A^{-1}$ just by taking $0 \leq \alpha_2 - \alpha_1 \leq \frac{1}{\lambda\,\alpha_1}$. The process can be repeated making $\lambda\,\alpha_n^n = \theta_n$ (mod 1) by perturbing the value of $\alpha_{n-1}$ less than $\lambda^{-1}\,\alpha_{n-1}^{1-n} < \lambda^{-1}\,A^{1-n}$, in such a way that the sequence $\alpha_n$ will converge to a limit $\alpha$ such that

$$\|\lambda\,\alpha^n - \theta_n\| \;\leq\; \sum_{k=1}^{\infty} A^{-k} \;=\; \frac{1}{A-1} \tag{2.6}$$

for $n = 1, 2, \ldots$ (actually, a little trick at the end of the proof allows us to halve the bound).

From here, the existence of uncountably many real numbers, so uncountably many transcendental numbers, whose powers are not u.d. mod 1, follows easily (corollary 3.6).

In Pisot's thesis ([8]) there are some general results related to the previous one. In particular, the application of his theorem II, on p. 215 of [8], would show that there are arbitrarily large numbers $\alpha$ such that:

$$\limsup_{n\to\infty} \|\alpha^n - \theta_n\| \;\leq\; \frac{1}{2\,(\alpha-1)} \tag{2.7}$$

Also, his theorem on p. 225 of [8] implies that for every $\varepsilon > 0$ there is some arbitrarily large $\alpha$ such that:

$$\|\alpha^n - \theta_n\| \; < \; \frac{1+\varepsilon}{2\,(\alpha-1)} \tag{2.8}$$

for $n \geq n_0(\varepsilon)$.

The next question that arises is if the sequence $\langle\alpha^n\rangle$ = fractional part of $\alpha^n$ uniquely determines $\alpha$. This can be stated as if $\alpha^n - \beta^n \in \mathbb{Z}$ for every $n \geq n_0$, for some $n_0$, implies that $\alpha = \beta$. The answer is yes, except in the trivial case when $\alpha$ and $\beta$ are both integers.

The techniques used to solve this problem are easy to generalize to the study of combinations $S_n = \lambda\,\alpha^n + \mu\,\beta^n$ (with $\lambda, \mu, \alpha, \beta \neq 0$ and $\alpha \neq \beta$) such that $S_n \in \mathbb{Z}$ for every $n \geq n_0$ for some $n_0$. The result now is that $\alpha$ and $\beta$ are rational integers or conjugate algebraic integers of degree two. Furthermore, $\lambda \in \mathbb{Q}(\alpha)$ and $\mu \in \mathbb{Q}(\beta)$. In general we do not have that $\mathbb{Q}(\lambda) = \mathbb{Q}(\alpha)$ and $\mathbb{Q}(\mu) = \mathbb{Q}(\beta)$. Of course, if $\alpha$ is a quadratic integer and $\beta = \bar{\alpha}$ denotes its conjugate, then $\alpha^n + \beta^n \in \mathbb{Z}$, whilst $\lambda = 1$ does not generate $\mathbb{Q}(\alpha)$. However, the equality does hold whenever $\lambda \neq \mu$. Pisot's thesis ([8]) also contains similar results, applicable to sequences of the form

$$S_n \; = \; \lambda_1(n)\,\alpha_1^n + \lambda_2(n)\,\alpha_2^n + \cdots + \lambda_m(n)\,\alpha_m^n \tag{2.9}$$

where $\alpha_i \in \mathbb{C}$ and the coefficients are polynomials with complex coefficients.

A generalization of those results to local fields is given in theorem 4.2. Note that $S_n \in \mathbb{Z}$ can be expressed as $S_n \in \mathbb{Q}$ and $|S_n|_v \leq 1$ for every non-Archimedean place $v$ in $\mathbb{Q}$. The generalization given in theorem 4.2 allows us to study what happens at each individual non-Archimedean place $v$. In that theorem, $\mathbb{Q}$ is replaced by some subfield $k'$ of a finite extension of a local field $k_v$, where $k_v$ = completion of a number field $k$ at a non-Archimedean place $v$, the $\alpha_i$'s and the coefficients of the $\lambda_i$'s are in an arbitrary extension of $k'$, and for $n \geq n_0$, $S_n$ is assumed to be in $k'$ and to have bounded $v$-adic absolute value. Under those hypothesis, the main conclusion is that $\alpha_i \in \overline{k'}$ and $|\alpha_i|_v \leq 1$ for every $i$, which is the local version of being algebraic integers. The proof runs along the following lines:

Since $|S_n|_v$ is bounded for $n \geq n_0$, it is possible to multiply it by some $b \in k'$ so that $c_n = b\,S_n$ verifies that $|c_n|_v \leq 1$, i.e.,

$$c_n \; \in \; \mathcal{O}_{k',v} \; = \; \{\, x \in k' \, : \, |x|_v \leq 1 \,\} \tag{2.10}$$

for $n \geq n_0$. Furthermore, it is well known that if

$$f(X) \; = \; \prod_{i=1}^{m} (X - \alpha_i)^{d_i} = X^M - r_1\,X^{M-1} - r_2\,X^{M-2} - \cdots - r_M \tag{2.11}$$

where $d_i \; = \; 1 + \deg \lambda_i$, $i \; = \; 1, 2, \ldots, m$, and $M \; = \; \sum_{i=1}^{m} d_i$, then $c_n$ verifies a recurrence relation of the form

$$c_n \; = \; r_1\,c_{n-1} + r_2\,c_{n-2} + \cdots + r_M\,c_{n-M} \tag{2.12}$$

Next, define the formal power series

$$F_t(X) \;=\; c_t + c_{t+1}\,X + c_{t+2}\,X^2 + c_{t+3}\,X^3 + \dots \tag{2.13}$$

which is in $\mathcal{O}_{k',v}[[X]]$. We have $F_t(X) = p_t(X)/q(X)$, where $p_t(X), q(X) \in k'[X]$ and $q(X) = X^M f(1/X)$ is the reciprocal of $f(X)$. Using the generalized Fatou's lemma (lemma 4.3), we conclude that $q(X) \in \mathcal{O}_{k',v}[X]$, so $f(X) \in \mathcal{O}_{k',v}[X]$. Since $f(X)$ is monic, we get that its roots verify $|\alpha_i|_v \le 1$.

Theorem 4.2 also shows that the coefficients of each $\lambda_i(X)$ are in $k'(\alpha_i)$. Even more, any $k'$-automorphism $\sigma$ in $k'(\alpha_1, \dots, \alpha_m)$, which acts as a permutation over the $\alpha_i$'s, induces the same permutation over the $\lambda_i$'s.

Finally, proposition 4.6 shows that the condition $|S-n|_v = O(1)$ can be weakened to a subexponential growth condition of the form $|S-n|_v = O(A^n)$ for every $A > 1$.

In the next sections we give details about the above results.

## 3. Distribution of powers modulo 1

**Definition 3.1.** For any real number $x$, we define:

(1) Integer part of $x$: $[x] = \max\{\, n \in \mathbb{Z} \,:\, n \le x \,\}$.

(2) Fractional part of $x$: $\langle x\rangle = x - [x]$.

(3) Nearest integer to $x$: $E(x) = \max\{\, n \in \mathbb{Z} \,:\, n \le x + 1/2 \,\}$.

(4) Residue of $x$ modulo 1: $\varepsilon(x) = x - E(x)$.

(5) Distance from $x$ the the nearest integer:

$$\|x\| = |\varepsilon(x)| = \min\{\, |x-n| \,:\, n \in \mathbb{Z} \,\}$$

**Definition 3.2.** A sequence of real numbers $\{x_n\}_{n=1}^{\infty}$ is said to be uniformly distributed modulo 1, abbreviated u.d. mod 1, if

$$\lim_{N\to\infty} \frac{1}{N} \sum_{n=1}^{N} \chi_{s,t}(x_n) \;=\; t - s \tag{3.1}$$

whenever $s < t < s+1$. Here:

$$\chi_{s,t}(x) \;=\; \begin{cases} 1 & \text{if} \quad s < x-n < t \quad \text{for some} \quad n \in \mathbb{Z} \\ 1/2 & \text{if} \quad s - x \in \mathbb{Z} \quad \text{or} \quad t - x \in \mathbb{Z} \\ 0 & \text{otherwise} \end{cases} \tag{3.2}$$

**Theorem 3.3** (Weyl criterion). *The sequence $\{x_n\}_{n=1}^{\infty}$ is u.d. mod 1 if and only if for every $h \in \mathbb{Z} \setminus \{0\}$:*

$$\lim_{N\to\infty} \frac{1}{N} \sum_{n=1}^{N} \exp(2\pi i\, h\, x_n) \;=\; 0 \tag{3.3}$$

*Proof.* See theorem 2.1 in [5]. □

**Theorem 3.4** (Weyl-Koksma's Metric Theorem).

*(a) Let $\alpha > 1$ be a real number; the sequence $\{\lambda\,\alpha^n\}_{n=1}^{\infty}$ is uniformly distributed modulo 1 for almost all real $\lambda$ (Weyl).*

*(b) Let $\lambda$ be a non zero real number; the sequence $\{\lambda\,\alpha^n\}_{n=1}^{\infty}$ is uniformly distributed modulo 1 for almost all real $\alpha > 1$ (Koksma).*

*Proof.* See [14] and [4]. See also [1], p. 71. □

**Theorem 3.5.** *Let $\{\theta_n\}$ be any sequence of real numbers. Then:*

(1) *Given any $\lambda \neq 0$ and $A > 1$, there exists an $\alpha$ such that:*

$$A \;\leq\; \alpha \;\leq A + \frac{A}{|\lambda|\,(A-1)} \tag{3.4}$$

*and for every $n \geq 1$*

$$\|\lambda\,\alpha^n - \theta_n\| \;\leq\; \frac{1}{2\,(A-1)} \tag{3.5}$$

(2) *Given any $\alpha > 1$ and $L \neq 0$, there exists a $\lambda$ (with the same sign as $L$) such that:*

$$|L| \;\leq\; |\lambda| \;\leq\; |L| + \frac{1}{\alpha - 1} \tag{3.6}$$

*and for every $n \geq 1$*

$$\|\lambda\,\alpha^n - \theta_n\| \;\leq\; \frac{1}{2\,(\alpha-1)} \tag{3.7}$$

*Proof.*

(1) There is no loss of generality in assuming that $\lambda > 0$. Otherwise, make $\lambda' = -\lambda$ and $\theta'_n = -\theta_n$, and use that $\|-x\| = \|x\|$.

Construct an increasing sequence $A = \alpha_0 \leq \alpha_1 \leq \alpha_2 \ldots$, by the following recursive rule ($n \geq 0$):

$$\begin{aligned} \alpha_{n+1} &= \lambda^{-1/(n+1)}\,(\lambda\,\alpha_n^{n+1} + \langle\theta_{n+1} - \lambda\,\alpha_n^{n+1}\rangle)^{1/(n+1)} \\ &= \lambda^{-1/(n+1)}\,(\theta_{n+1} - [\theta_{n+1} - \lambda\,\alpha_n^{n+1}])^{1/(n+1)} \end{aligned} \tag{3.8}$$

where $[x]$ = integer part of $x$, and $\langle x\rangle = x - [x]$.

Since $\langle\theta_{n+1} - \lambda\,\alpha_n^{n+1}\rangle \geq 0$ we have that $\alpha_{n+1} \geq \alpha_n$, so the sequence is non decreasing. Also, by construction $\langle\lambda\,\alpha_n^n\rangle = \langle\theta_n\rangle$ for every $n \geq 1$.

Now, for every $1 \leq m \leq n$:

(3.9)

$$\begin{aligned}
\alpha_{n+1}^m - \alpha_n^m &= \lambda^{-m/(n+1)} \, (\lambda \, \alpha_n^{n+1} + \langle \theta_{n+1} - \lambda \, \alpha_n^{n+1} \rangle)^{m/(n+1)} - \alpha_n^m \\
&\leq \lambda^{-m/(n+1)} \, \Big( (\lambda \, \alpha_n^{n+1} + 1)^{m/(n+1)} - (\lambda \, \alpha_n^{n+1})^{m/(n+1)} \Big) \\
&= \lambda^{-m/(n+1)} \, \frac{m}{n+1} \int_0^1 (\lambda \, \alpha_n^{n+1} + x)^{(m-n-1)/(n+1)} \, dx \\
&\leq \lambda^{-m/(n+1)} \, \frac{m}{n+1} \, (\lambda \, \alpha_n^{n+1})^{(m-n-1)/(n+1)} \\
&< \lambda^{-1} \, \alpha_n^{m-n-1} \; \leq \; \lambda^{-1} \, A^{m-n-1}
\end{aligned}$$

For $m = 1$ we have $\alpha_{n+1} - \alpha_n < \lambda^{-1} \, A^{-n}$, thus:

$$\begin{aligned}
(3.10) \quad \alpha_N - \alpha_0 \; &= \; \sum_{n=0}^{N-1} (\alpha_{n+1} - \alpha_n) \\
&= \; \sum_{n=0}^{N-1} \lambda^{-1} \, A^{-n} \; < \; \sum_{n=0}^{\infty} \lambda^{-1} \, A^{-n} \; = \; \frac{A}{\lambda \, (A-1)}
\end{aligned}$$

So, for every $N \geq 0$ we have $A \leq \alpha_N < A + \frac{A}{\lambda \, (A-1)}$. Hence, $\alpha = \lim_{N \to \infty} \alpha_N$ exists, and

$$A \; \leq \; \alpha \; \leq \; A + \frac{A}{\lambda \, (A-1)}$$

which is (3.4).

Now, for $N > m \geq 1$:

$$\begin{aligned}
(3.11) \quad 0 \; \leq \lambda \, \alpha_N^m - \lambda \, \alpha_m^m \; &= \; \sum_{n=m}^{N-1} \lambda \, (\alpha_{n+1}^m - \alpha_n^m) \\
&\leq \; \sum_{n=m}^{N-1} A^{m-n-1} \; < \; \sum_{n=m}^{\infty} A^{m-n-1} \; = \; \frac{1}{A-1}
\end{aligned}$$

Now let $N \to \infty$:

$$(3.12) \qquad 0 \; \leq \lambda \, \alpha^m - \lambda \, \alpha_m^m \; \leq \; \frac{1}{A-1}$$

Since $\langle \lambda \, \alpha_m^m \rangle = \langle \theta_m \rangle$, we get:

$$(3.13) \qquad \langle \lambda \, \alpha^m - \theta_m \rangle \; \leq \; \frac{1}{A-1}$$

To get the desired result, apply the previous result to the sequence $\theta'_n = \theta_n - \frac{1}{2 \, (A-1)}$, i.e., for $n \geq 1$:

$$(3.14) \qquad 0 \; \leq \; \langle \lambda \, \alpha^n - \theta_n + \frac{1}{2 \, (A-1)} \rangle \; \leq \; \frac{1}{A-1}$$

Subtracting $\frac{1}{2\,(A-1)}$ we get:

$$\left|\langle\lambda\,\alpha^n - \theta_n + \frac{1}{2\,(A-1)}\rangle - \frac{1}{2\,(A-1)}\right| \;\leq\; \frac{1}{2\,(A-1)} \tag{3.15}$$

or, using $\langle x\rangle = x - [x]$:

$$\left|\lambda\,\alpha^n - \theta_n - \left[\lambda\,\alpha^n - \theta_n + \frac{1}{2\,(A-1)}\right]\right| \;\leq\; \frac{1}{2\,(A-1)} \tag{3.16}$$

Hence:

$$\|\lambda\,\alpha^n - \theta_n\| \;\leq\; \frac{1}{2\,(A-1)}$$

which is (3.5).

(2) There is no loss of generality in assuming that $L > 0$. Otherwise, make $L' = -L$, $\lambda' = -\lambda$ and $\theta'_n = -\theta_n$, and use that $\|-x\| = \|x\|$.

Construct an increasing sequence $L = \lambda_0 \leq \lambda_1 \leq \lambda_2 \leq \dots$, by the following recursive rule for $n \geq 0$:

$$\begin{aligned} \lambda_{n+1} &= \lambda_n + \alpha^{-n-1}\,\langle\theta_{n+1} - \lambda_n\,\alpha^{n+1}\rangle \\ &= \alpha^{-n-1}\,(\theta_{n+1} - [\theta_{n+1} - \lambda\,\alpha^{n+1}]) \end{aligned} \tag{3.17}$$

Since $\langle\theta_{n+1} - \lambda_n\,\alpha^{n+1}\rangle \geq 0$, we have that $\lambda_{n+1} \geq \lambda_n$ $(n \geq 0)$, so the sequence is actually increasing. Also, by construction $\langle\lambda_n\,\alpha^n\rangle = \langle\theta_n\rangle$ for every $n \geq 1$. Furthermore, since $\langle\theta_{n+1} - \lambda_n\,\alpha^{n+1}\rangle < 1$:

$$0 \;\leq\; \lambda_{n+1} - \lambda_n \;<\; \alpha^{-n-1} \tag{3.18}$$

so, for $N \geq 1$:

$$\begin{aligned} 0 \;\leq\; \lambda_N - \lambda_0 \;&=\; \sum_{n=0}^{N-1}(\lambda_{n+1} - \lambda_n) \\ &<\; \sum_{n=0}^{N-1}\alpha^{-n-1} \;<\; \sum_{n=0}^{\infty}\alpha^{-n-1} \;=\; \frac{1}{\alpha-1} \end{aligned} \tag{3.19}$$

Hence, $L \leq \lambda_N < L + \frac{1}{\alpha-1}$, so $\lambda = \lim_{N\to\infty}\lambda_N$ exists, and:

$$L \;\leq\; \lambda \;\leq L + \frac{1}{\alpha-1}$$

which is (3.6).

Now, for $N > m \geq 1$:

$$\begin{aligned} 0 \;\leq\; \lambda_N\,\alpha^m - \lambda_m\,\alpha^m \;&=\; \alpha^m\sum_{n=m}^{N-1}(\lambda_{n+1} - \lambda_n) \\ &<\; \alpha^m\sum_{n=m}^{N-1}\alpha^{-n-1} \;<\; \sum_{n=m}^{\infty}\alpha^{m-n-1} \;=\; \frac{1}{\alpha-1} \end{aligned} \tag{3.20}$$

Letting $N \to \infty$:

$$0 \;\leq\; \lambda\,\alpha^m - \lambda_m\,\alpha^m \;\leq\; \frac{1}{\alpha-1} \tag{3.21}$$

Since $\langle \lambda_m \alpha^m \rangle = \langle \theta_m \rangle$, we get:

$$\langle \lambda \alpha^m - \theta_m \rangle \leq \frac{1}{\alpha - 1} \tag{3.22}$$

To get the desired result, apply the previous result to the sequence $\theta'_n = \theta_n - \frac{1}{\alpha-1}$, and apply the same reasoning as in part (1).

□

**Corollary 3.6.** *For any $\lambda \neq 0$ there are uncountably many numbers $\alpha > 1$ such that $\lambda \alpha^n$ is not uniformly distributed modulo 1. Hence, there are transcendental numbers with that property.*

*Proof.* For every subset $I \subseteq \mathbb{Z}^+$ take an $\alpha_I > 1$ such that $\|\lambda \alpha_I^n - 1/4\| \leq 1/5$ if $n \in I$ and $\|\lambda \alpha_I^n - 3/4\| \leq 1/5$ if $n \notin I$. Such $\alpha_I$ exists by theorem 3.5. We have that $I \neq I' \implies \alpha_I \neq \alpha_{I'}$, so the map $I \mapsto \alpha_I$ is injective. Since there are uncountably many subsets in $\mathbb{Z}^+$, the result follows. □

## 4. Combinations of powers

*Remark 4.1* (Notation). In this section, $k$ will be a number field, $v$ a non-trivial, non-Archimedean place in $k$, $k_v$ the completion of $k$ at $v$, $\overline{k}_v$ an algebraic closure of $k_v$, $\Omega_v$ the completion of $\overline{k}_v$, $k' \subseteq k_v(\gamma) \subseteq \overline{k}_v$ any subfield of a finite extension $k_v(\gamma)$ of $k_v$ ($\gamma \in \overline{k}_v$), $\overline{k'} \subseteq \overline{k}_v$ the set of elements of $\overline{k}_v$ that are algebraic over $k'$, $K$ any field extension of $\overline{k'}$, and $\mathcal{O}_{k',v}$ the ring $\mathcal{O}_{k',v} = \{ x \in k' : |x|_v \leq 1 \}$.

**Theorem 4.2.** *Let $\lambda_1(X), \lambda_2(X), \ldots, \lambda_m(X)$ be polynomials in $K[X] \setminus \{0\}$. Assume that $\alpha_1, \alpha_2, \ldots, \alpha_m$ are distinct elements from $K \setminus \{0\}$. Let $S_n$ be the following sum:*

$$S_n = \sum_{i=1}^{m} \lambda_i(n)\, \alpha_i^n \tag{4.1}$$

*Also assume that there is an integer $n_0$ such that for every $n \geq n_0$, $S_n \in k'$. Then for every $n \in \mathbb{Z}$, $S_n \in k'$, and:*

*(i) $\alpha_1, \alpha_2, \ldots, \alpha_m$ are the roots of a polynomial $f(X) \in k'[X]$, so they are in $\overline{k'}$. Moreover, if $|S_n|_v = O(1)$, then $|\alpha_i|_v \leq 1$, $i = 1, 2, \ldots, m$.*

*(ii) For every $i = 1, 2, \ldots, m$: $k'(\lambda_i) \subseteq k'(\alpha_i)$, where $k'(\lambda_i)$ is the field generated by the coefficients of $\lambda_i(X)$.*

*(iii) For every automorphism $\sigma \in \mathrm{Gal}(k'(\alpha_1, \ldots, \alpha_m)/k')$ such that $\sigma(\alpha_i) = \alpha_{\pi(i)}$, $i = 1, \ldots, m$, where $\pi$ is a permutation of $\{1, \ldots, m\}$, we have that $\sigma(\lambda_i(X)) = \lambda_{\pi(i)}(X)$, where by definition $\sigma(\sum_i a_i X^i) = \sum_i \sigma(a_i) X^i$.*

*(iv) If $\lambda_1(X), \ldots, \lambda_m(X)$ are distinct polynomials, then $k'(\alpha_i) = k'(\lambda_i)$ for every $i$.*

The proof of this theorem requires several lemmas.

**Lemma 4.3** (Generalized Fatou's Lemma). *Let $A$ be a Dedekind ring and $F$ a rational series in $A[[X]]$, i.e., $F = p/q$ for some $p, q \in A[X]$. Then there exist two polynomials $P, Q \in A[X]$ such that $F = P/Q$, where $P$ and $Q$ are relatively prime and $Q(0) = 1$.*

*Proof.* See [1], p. 15, theorem 1.3. □

**Lemma 4.4.** *Let $\{c_n\}_{n=-\infty}^{\infty}$ a set of elements from $K$ such that $c_n \in k'$ for every $n \geq n_0$, and verifying the following recurrence relation of order M:*

$$c_n \;=\; r_1\, c_{n-1} + r_2\, c_{n-2} + \cdots + r_M\, c_{n-M} \tag{4.2}$$

*for every $n \in \mathbb{Z}$, where $r_1, r_2, \ldots, r_M$ are in $K$, $r_M \neq 0$. Then:*

*(i) The coefficients $r_1, r_2, \ldots, r_M$ are in $k'$, and for every $n \in \mathbb{Z}$, $c_n \in k'$.*

*(ii) If $c_n \in \mathcal{O}_{k',v}$ for every $n \geq n_0$, then the coefficients $r_1, r_2, \ldots, r_M$ are all in $\mathcal{O}_{k',v}$.*

*Proof.*

(i) Let $C_n$ and $R$ be the matrices:

$$C_n \;=\; \begin{pmatrix} c_n & c_{n+1} & \cdots & c_{n+M-1} \\ c_{n+1} & c_{n+2} & \cdots & c_{n+M} \\ \vdots & \vdots & \ddots & \vdots \\ c_{n+M-1} & c_{n+M} & \cdots & c_{n+2M-2} \end{pmatrix} \tag{4.3}$$

and

$$R \;=\; \begin{pmatrix} 0 & 1 & 0 & \ldots & 0 \\ 0 & 0 & 1 & \ldots & 0 \\ \vdots & \vdots & \vdots & \ddots & \vdots \\ 0 & 0 & 0 & \ldots & 1 \\ r_M & r_{M-1} & r_{M-2} & \ldots & r_1 \end{pmatrix} \tag{4.4}$$

We have that $C_{n+1} = R\, C_n$. Since the recurrence relation is of order M, $C_n$ is non singular. On the other hand, $R = C_{n+1}\, C_n^{-1}$. Since the elements of $C_n$ are in $k'$ for $n \geq n_0$, the entries of $R$, and those of $R^{-1}$, will be in $k'$. Since $C_{n-1} = R^{-1}\, C_n$, we get that the entries of $C_n$ will be in $k'$ also for $n < n_0$.

(ii) For each $t \geq n_0$ define the formal power series

$$F_t(X) \;=\; \sum_{n=0}^{\infty} c_{t+n}\, X^n \tag{4.5}$$

which is in $\mathcal{O}_{k',v}[[X]]$. We have $F_t(X) = p_t(X)/q(X)$, where $p_t(X), q(X) \in k'[X]$ are the following:

$$p_t(X) \;=\; \sum_{j=0}^{M-1} \Big( c_{t+j} - \sum_{i=1}^{j} r_i\, c_{t+j-i} \Big)\, X^j \tag{4.6}$$

$$q(X) \;=\; 1 - r_1\, X - r_2\, X^2 - \cdots - r_M\, X^M \tag{4.7}$$

This can be checked by multiplying $F_t(X)$ by $q_t(X)$ and using the recurrence relation, which gives $F_t(X)\, q(X) = p_t(X)$ (see [10] and [11]).

Now we will prove that $p_t(X)$ and $q(X)$ are relatively prime. To do so, we will see that they cannot have any common root (in $\overline{k'}$). In fact, assume that $\alpha$ is a common root of $p_{t_0}(X)$ and $q(X)$ for some $t_0 \geq n_0$, i.e.: $p_{t_0}(\alpha) = q(\alpha) = 0$. Since $q(0) = 1$, then $\alpha \neq 0$. Now we have:

$$X\, F_{t_0+1}(X) = F_{t_0}(X) - c_{t_0} \tag{4.8}$$

so:

$$\begin{aligned} X\, p_{t_0+1}(X) &= X\, q(X)\, F_{t_0+1}(X) \\ &= q(X)\, (F_{t_0}(X) - c_{t_0}) = p_{t_0}(X) - c_{t_0}\, q(X) \end{aligned} \tag{4.9}$$

Hence $p_{t_0+1}(\alpha) = 0$, which means that $\alpha$ is also a root of $p_{t_0+1}(X)$. By induction we get that $p_t(\alpha) = 0$ for every $t \geq t_0$. Grouping the terms of $p_t(X)$ with respect to $c_t, c_{t+1}, \dots, c_{t+M-1}$, we get:

$$p_t(X) = \sum_{j=0}^{M-1} a_j(X)\, c_{t+j} \tag{4.10}$$

where

$$a_j(X) = X^j \Big( 1 - \sum_{i=1}^{M-j-1} r_i\, X^i \Big) \tag{4.11}$$

Note that $a_0(X), a_1(X), \dots, a_{M-1}(X)$ do not depend on t. On the other hand $p_t(\alpha) = 0$ implies

$$\sum_{j=0}^{M-1} a_j(\alpha)\, c_{t+j} = 0 \tag{4.12}$$

for every $t \geq t_0$. Note that $a_{M-1}(\alpha) = \alpha^{M-1} \neq 0$, so $a_0(\alpha), a_1(\alpha), \dots, a_{M-1}(\alpha)$ are not all zero, and (4.12) means that the columns of the matrix $C_{t_0}$ are linearly dependent, so $\det C_{t_0} = 0$, which contradicts the fact that $C_{t_0}$ is non singular. Hence, the hypothesis that $p_t(X)$ and $q(X)$ have a common root has to be false. This proves that $p_t(X)$ and $q(X)$ are relatively prime.

By (generalized Fatou's) lemma 4.3, and taking into account that $\mathcal{O}_{k',v}$ is a Dedekind ring (see, for instance, [7], chap. 5), we get that there exist two relatively prime polynomials $P_t(X)$ and $Q_t(X)$ in $\mathcal{O}_{k',v}[X]$ such that $F_t(X) = P_t(X)/Q_t(X)$ and $Q_t(0) = 1$. Hence: $p_t(X)\, Q_t(X) = q(X)\, P_t(X)$. By unique factorization of polynomials in $k'[X]$, there is a $u \in k'$ such that $P_t(X) = u\, p_t(X)$ and $Q_t(X) =$ $u\,q_t(X)$. Since $Q_t(0) = q(0) = 1$, we get that $u = 1$, so $P_t(X) = p_t(X)$ and $Q_t(X) = q(X)$. Hence, the coefficients of $q(X)$ are in $\mathcal{O}_{k',v}$.

□

*Proof of Theorem 4.2.* Let $f(X) \in K[X]$ be the polynomial

(4.13)
$$f(X) = \prod_{i=1}^{m} (X - \alpha_i)^{d_i} = X^M - r_1\, X^{M-1} - r_2\, X^{M-2} - \cdots - r_M$$

where $d_i = 1 + \deg \lambda_i$, $i = 1, 2, \ldots, m$, and $M = \sum_{i=1}^{m} d_i$. Note that this is just the reciprocal of the polynomial $q(X)$ defined in the proof lemma 4.4. It is well known (see [10] and [11]) that $S_n$ verifies a recurrence relation of order $M$, of the form:

$$S_n \;=\; r_1\, S_{n-1} + r_2\, S_{n-2} + \cdots + r_M\, S_{n-M} \tag{4.14}$$

By part (i) of lemma 4.4 we get that $S_n \in k'$ for every $n \in \mathbb{Z}$.

Now we prove the remaining results.

(i) Since $\alpha_1, \ldots, \alpha_m$ are the roots of $f(X) \in k'[X]$, it is clear that they are in $\overline{k'}$.

Assume that $|S_n|_v = O(1)$. We have that there are integers $B > 0$ and $n_0$ such that for every $n \geq n_0$, $|S_n|_v \leq B$. Let $b \in k'$ be any element such that $|b|_v \leq 1/B$. Then for every $n \geq n_0$, $|b\, S_n|_v \leq 1$. Putting $c_n = b\, S_n$ we get that $c_n \in \mathcal{O}_{k',v}$ for every $n \geq n_0$, and also verifies the recurrence:

$$c_n \;=\; r_1\, c_{n-1} + r_2\, c_{n-2} + \cdots + r_M\, c_{n-M} \tag{4.15}$$

By part (ii) of lemma 4.4 we get that $r_1, r_2, \ldots, r_M$ are in $\mathcal{O}_{k',v}$, i.e., $|r_i|_v \leq 1$ for $i = 1, 2, \ldots, M$. We already know that $\alpha_i \in \overline{k'}$, and using the well known equality:

$$|r_0|_v \prod_{l=1}^{m} \max\{1, |\alpha_l|_v\}^{d_l} \;=\; \max\{\, |r_l|_v \,:\, 0 \leq l \leq M \,\} \tag{4.16}$$

($r_0 = 1$), we get that actually $|\alpha_i|_v \leq 1$.

(ii) For $i = 1, 2, \ldots, m$, we have $\lambda_i(X) \;=\; \sum_{j=0}^{d_j - 1} a_{ij}\, X^j$, where the coefficients $a_{ij}$ are in $K$. If $d = \max\{d_1, d_2, \ldots, d_m\}$, then we can write $\lambda_i(X) \;=\; \sum_{j=0}^{d-1} a_{ij}\, X^j$, where $a_{ij} = 0$ for $d_i \leq j < d$. Now consider the following matrices:

$$\vec{S} \;=\; \begin{pmatrix} S_0 \\ S_1 \\ \vdots \\ S_{md-1} \end{pmatrix}, \qquad \vec{L} \;=\; \begin{pmatrix} \vec{L}_1 \\ \vec{L}_2 \\ \vdots \\ \vec{L}_m \end{pmatrix} \tag{4.17}$$

where the $\vec{L}_i$'s are blocks of the form:

$$\vec{L}_i \;=\; \begin{pmatrix} a_{i\,0} \\ a_{i\,1} \\ \vdots \\ a_{i\,d-1} \end{pmatrix} \tag{4.18}$$

and:

$$A = (A_1 \, A_2 \, \dots \, A_m) \tag{4.19}$$

where the $A_i$'s are blocks of the form:

$$A_i = \begin{pmatrix} 1 & 0 & \dots & 0 \\ \alpha_i & \alpha_i & \dots & \alpha_i \\ {\alpha_i}^2 & 2{\alpha_i}^2 & \dots & 2^{d-1}\,{\alpha_i}^2 \\ {\alpha_i}^3 & 3{\alpha_i}^3 & \dots & 3^{d-1}\,{\alpha_i}^3 \\ \vdots & \vdots & \ddots & \vdots \\ {\alpha_i}^{(md-1)} & (md-1)\,{\alpha_i}^{(md-1)} & \dots & (md-1)^{(d-1)}\,{\alpha_i}^{(md-1)} \end{pmatrix} \tag{4.20}$$

i.e., the element of $A_i$ in row $n$, column $j$, is $n^j\,\alpha_i^n$, for $0 \leq n \leq md-1$ and $0 \leq j \leq d-1$ (by convention, $0^0 = 1$). Note that $A$ is a square matrix.

It is easy to check that $\vec{S} = A\,\vec{L}$. Furthermore, it is known that $\det A \neq 0$ (see lemma 8.5.1 in [13], 177-182). Next, considering the coefficients $a_{ij}$ as unknowns of a system of $md$ linear equations and $md$ unknowns, and using Cramer's rule, we get

$$a_{ij} = \det(A'_{ij}) / \det(A) \tag{4.21}$$

where $A'_{ij}$ is matrix $A$ with the $j$-th column of block $A_i$ substituted by $\vec{S}$. Considering $a_{ij}$ as a function of $\alpha_1, \alpha_2, \dots, \alpha_{i-1}, \alpha_{i+1}, \dots, \alpha_m$, we note that interchanging any two of them, say $\alpha_l$ and $\alpha_r$ $(l, r \neq i)$, does not change the value of $a_{ij}$, so $a_{ij}$ is a symmetric rational function of the $\alpha_l$'s $(l \neq 1)$, with coefficients in $k'(\alpha_i)$. Hence, $a_{ij}$ is a rational function in

$$k'(\alpha_i)(\phi_{i1}, \phi_{i2}, \dots, \phi_{i,i-1}, \phi_{i,i+1}, \dots, \phi_{i,m}) \tag{4.22}$$

where $\phi_{il}$ for $1 \leq l \leq m$, $l \neq i$, are the elementary symmetric functions of $\{\alpha_l\}_{l \neq i}$. In other words, the $\phi_{il}$'s are $\pm$ the coefficients of the polynomial

$$\prod_{\substack{l=1 \\ l \neq i}}^{m} (X - \alpha_l) \;=\; \frac{\prod_{l=1}^{m}(X - \alpha_l)}{X - \alpha_i} \tag{4.23}$$

Since $f(X) = \prod_{l=1}^{m}(X - \alpha_l)^{d_l} \in k'[X]$, and in characteristic zero every irreducible polynomial is separable, it is clear that the polynomial in the numerator of (4.23) is in $k'[X]$. Hence, the quotient is in $k'(\alpha_i)[[X]]$, and actually in $k'(\alpha_i)[X]$, since it equals a polynomial. This implies that the $\phi_{il}$'s are in $k'(\alpha_i)$, hence $a_{ij} \in k'(\alpha_i)$. This proves the desired result.

(iii) We use again the expression $a_{ij} = \det(A'_{ij}) / \det(A)$. If $\sigma(\alpha_i) = \alpha_{\pi(i)}$, then:

$$\sigma(\det(A_1\, A_2\, \dots\, A_m)) = \det(A_{\pi(1)}\, A_{\pi(2)}\, \dots\, A_{\pi(m)}) \tag{4.24}$$

The effect of $\sigma$ on $\det(A'_{ij})$ is similar, but now $S$ will be in the $j$-th column of block $A_{\pi(i)}$. Hence, $\sigma(\det(A'_{ij}))/\sigma(\det(A)) = a_{\pi(i),j}$, i.e., $\sigma(a_{ij}) = a_{\pi(i),j}$, which proves the result.

(iv) From (ii) we have that $k'(\lambda_i) \subseteq k'(\alpha_i)$, so we will prove the other containment.

If $\lambda_1(X), \dots, \lambda_m(X)$ are distinct non zero polynomials, then for some $n \in \mathbb{Z}$ the numbers $\lambda_1(n), \dots, \lambda_m(n)$ are distinct and non zero. For $j \in \mathbb{Z}$, consider the following sums:

$$T_j = \sum_{i=1}^{m} \lambda_i(n)^j \, \alpha_i \tag{4.25}$$

It is easy to see that for every $j \in \mathbb{Z}$, $T_j \in k'$. In fact, take any automorphism $\sigma \in \mathrm{Gal}(k'(\alpha_1, \dots, \alpha_m)/k')$ such that $\sigma(\alpha_i) = \alpha_{\pi(i)}$ for $i = 1, 2, \dots, m$, where $\pi$ is some permutation of $\{1, 2, \dots, m\}$. Then:

$$\sigma(T_j) = \sum_{i=1}^{m} \lambda_{\pi(i)}(n)^j \, \alpha_{\pi(i)} = \sum_{i=1}^{m} \lambda_i(n)^j \, \alpha_i = T_j \tag{4.26}$$

So, $T_j$ is invariant for $\sigma$. This implies that $T_j \in k'$.

Now, consider each $\alpha_i$ as a polynomial of degree zero. Reasoning as in the proof of (ii) with the roles of $\lambda_i$ and $\alpha_i$ interchanged, we get that $k'(\alpha_i) \subseteq k'(\lambda_i(n))$. Since $\lambda_i(n) \in k'(\lambda_i)$, we get the desired result.

□

Next corollary is the global version of part (i) of theorem 4.2. Here, we take $k' = k$, and $\mathcal{O}_k = \bigcap_{\text{finite } v} \mathcal{O}_{k',v}$.

**Corollary 4.5.** *Under the same hypothesis as in theorem 4.2 with $k' = k$, if there is an element $b \in k$ such that $b\,S_n \in \mathcal{O}_k$ for every $n \geq n_0$, then $\alpha_i$ $(i = 1, 2, \dots, m)$ are algebraic integers.*

*Proof.* $b\,S_n \in \mathcal{O}_k$ for every $n \geq n_0$ implies that $S_n \in k$ and $|S_n|_v = O(1)$ at every finite place $v$, so that by theorem 4.2 we get $\alpha \in \overline{k}$ and $|\alpha_i|_v \leq 1$ for every finite $v$, hence, the $\alpha_i$'s are algebraic integers. □

Next proposition shows that the condition $|S_n|_v = O(1)$ in theorem 4.2 can be weakened to a growth condition.

**Proposition 4.6.** *Under the hypothesis of theorem 4.2, if $|S_n|_v = O(A^n)$ for every $A > 1$, then $|S_n|_v = O(1)$.*

*Proof.* Let $\pi$ be any element in $\overline{k'}$ such that $\pi^2$ is a prime generator of the ideal

$$\mathcal{M} = \{\, x \in k'(\alpha_1, \dots, \alpha_m) \,:\, |x|_v < 1 \,\} \tag{4.27}$$

Take $A = |\pi|_v^{-1}$. Let $\alpha_i' = \alpha_i\,\pi$, and let $S_n'$ be:

$$S_n' = \sum_{i=1}^{m} \lambda_i(n)\,(\alpha_i')^n = \sum_{i=1}^{m} \lambda_i(n)\,(\alpha_i\,\pi)^n = S_n\,\pi^n \tag{4.28}$$

By hypothesis, $|S'_n|_v = |S_n|_v\,|\pi|_v^n = |S_n|_v/A^n = O(1)$. Furthermore, we have that $S'_n \in k'(\pi) \subseteq k_v(\gamma,\pi)$. Since $k_v(\gamma,\pi)$ is a finite extension of $k_v$, theorem 4.2 applies, hence, by part (i) of that theorem, for every $i$, $|\alpha'_i|_v = |\alpha_i\,\pi|_v \le 1$, thus $|\alpha_i|_v \le |\pi|_v^{-1}$, which is strictly less than $|\pi^2|_v^{-1} = \min\{\,|x|_v^{-1} \,:\, x \in \mathcal{M}\setminus\{0\}\,\}$, so $|\alpha_i|_v \le 1$, and from here the result follows.

□

## 5. Other topics

**5.1. Asymptotic results.** We mentioned in section 1 that P.V. numbers belong to the exceptional set of Koksma's theorem, because if $\alpha$ is a P.V. number, then $\lim_{n\to\infty}\|\alpha^n\| = 0$ (geometrically). We could ask if there are other real numbers ($> 1$), besides the P.V. numbers, with this property. In other words, given a real number $\alpha > 1$, does $\lim_{n\to\infty}\|\alpha^n\| = 0$ imply that $\alpha$ is a P.V. number? The answer is not known so far, however it is possible to prove that $\alpha$ is actually a P.V. number when any of the following (increasingly weak) conditions is added ([1], section 5.4):

(1) $\alpha$ is algebraic, and there exists a real $\lambda \neq 0$ such that:

$$\lim_{n\to\infty}\|\lambda\,\alpha^n\| \;=\; 0 \tag{5.1}$$

(2) There exists a real $\lambda \neq 0$ such that:

$$\sum_{n=1}^{\infty}\|\lambda\,\alpha^n\|^2 \;<\; +\infty \tag{5.2}$$

(3) There exists a real $\lambda \neq 0$ such that:

$$\|\lambda\,\alpha^n\|^2 \;=\; o(n^{-1/2}) \tag{5.3}$$

(4) There exists a real $\lambda \neq 0$ such that:

$$\|\lambda\,\alpha^n\|^2 \;\le\; \frac{a}{\sqrt{n}} \quad (\forall n \ge n_0) \tag{5.4}$$

for some $0 < a < 1/2\sqrt{2}(\alpha+1)^2$ and some integer $n_0$.

The proofs rest on the rationality of the power series $f(X) = \sum_{n=0}^{\infty} u_n\,X^n$, where $u_n =$ nearest integer to $\lambda\,\alpha^n$. Next, by Fatou's lemma, $f(X) = A(X)/Q(X)$ for some relatively prime $A, Q \in \mathbb{Z}[X]$, such that $Q(0) = 1$. On the other hand we have that:

$$f(X) \;=\; \frac{A(X)}{Q(X)} \;=\; \frac{\lambda}{1-\alpha\,X} + \varepsilon(X) \tag{5.5}$$

where $\varepsilon(X) = \sum_{n=0}^{\infty}\varepsilon_n\,X^n$; $\varepsilon_n = u_n - \lambda\,\alpha^n \in [-1/2, 1/2)$. Since $\varepsilon_n \to 0$, we get that the meromorphic function $\varepsilon(z)$ has no pole on the close unit disk $|z| \le 1$, so $Q(z)$ has a single zero in the close unit disk. From here it is easy to see that $\alpha$ is a P.V. number ([1], sec. 5.4.; [12], 4-10).

**5.2. Generalizations to adeles.** Let $\mathbb{A}$ be the ring of adeles of $\mathbb{Q}$, and $I$ a finite set of places in $\mathbb{Q}$. The $I$-adele ring of $\mathbb{Q}$ is defined as:

$$\mathbb{A}_I \;=\; \{\, x \in \mathbb{A} \,:\, x_p = 0 \text{ for } p \notin I \,\} \tag{5.6}$$

Which is isomorphic to $\prod_{p\in I} \mathbb{Q}_p$ and contains a field canonically isomorphic to $\mathbb{Q}$ (which will also be designed $\mathbb{Q}$). Let $\mathbb{Q}^I$ be:

$$\mathbb{Q}^I \;=\; \{\, x \in \mathbb{Q} \,:\, |x|_p \le 1 \text{ for } p \notin I \cup \{\infty\} \,\} \tag{5.7}$$

Then we have that $\mathbb{Q}^I$ is a Dedekind ring, and:

(1) The field $\mathbb{Q}$ is dense in $\mathbb{A}_I$.

(2) $\mathbb{Q}^I$ is a discrete subring of $\mathbb{A}_I$, and the quotient $\mathbb{A}_I/\mathbb{Q}^I$ is locally compact.

Let $F_I$ be the set:

$$F_I \;=\; [-\frac{1}{2}, \frac{1}{2}) \times \prod_{p\in I\setminus\{\infty\}} \mathcal{O}_p \;\cong\; \mathbb{A}_I/\mathbb{Q}^I \tag{5.8}$$

then every element $x \in \mathbb{A}_I$ can be expressed in one and only one way as $x = E(x) + \varepsilon(x)$, with $E(x) \in \mathbb{Q}^I$ and $\varepsilon(x) \in F_I$ (Artin decomposition, see [1], theorem 10.1.2). Here, $E(x)$ plays the role of "nearest integer" to $x$, and $\varepsilon(x)$ is the residue of $x$ modulo $\mathbb{Q}^I$. Note that this theory becomes the usual theory of distribution modulo 1 when $I = \{\infty\}$.

In this setting we may define the concept of uniform distribution by using Weyl's criterion:

**Definition 5.1.** A sequence $\{x_n\}_{n=1}^{\infty}$ in $\mathbb{A}_I$ is uniformly distributed modulo $\mathbb{Q}^I$ if for all $a \in \mathbb{Q}^I$:

$$\lim_{N\to\infty} \frac{1}{N} \sum_{n=1}^{N} \exp(2\pi i\, \varepsilon_\infty(a\, x_n)) \;=\; 0 \tag{5.9}$$

Other possible definitions of uniform distribution over rings of adeles can be found in [3].

The version of Koksma's theorem here is the following ([1], theorem 10.1.6):

**Theorem 5.2.**

*(i) Let $\alpha \in \mathbb{A}_I$, with $|\alpha|_p > 1$ for every $p \in I$. Then the sequence $\{\lambda\, \alpha^n\}_{n=1}^{\infty}$ is uniformly distributed modulo $\mathbb{Q}^I$ for almost every invertable element $\lambda$ in $\mathbb{A}_I$.*

*(ii) Let $\lambda$ be an invertable element of $\mathbb{A}_I$. Then the sequence $\{\lambda\, \alpha^n\}_{n=1}^{\infty}$ is uniformly distributed modulo $\mathbb{Q}^I$ for almost all $\alpha \in \mathbb{A}_I$ with $|\alpha|_p > 1$ for every $p \in I$.*

**5.3. Generalization to fields of formal power series.** Let $k$ be any finite field, $\mathcal{Z} = k[X]$ and $\mathcal{F} = k(X)$. If $|k| = q$, then the following are absolute values on $\mathcal{F}$ ([1], chap. 12):

(1) If $f, g \in \mathcal{Z} \setminus \{0\}$, define $|f/g|_\infty = q^{(\deg f - \deg g)}$.

(2) If $v$ is a prime polynomial in $\mathcal{Z}$ and $f, g \in \mathcal{Z} \setminus \{0\}$ are relatively prime to $v$, then $|v^h\, f/g|_v = q^{-h}$.

Let $\mathcal{F}_v$ be the completion of $\mathcal{F}$ at place $v$, and $\mathcal{Z}_v = \{\, x \;:\; |x|_v \leq 1 \,\}$ the valuation ring of $\mathcal{F}_v$. Then we have that $\mathcal{F}_\infty = k\{X^{-1}\}$ = formal Laurent series of the form $\sum_{n=-\infty}^{h} a_n\, X^n$, $a_n \in k$. Furthermore, every element $x \in \mathcal{F}_\infty$ can be written in a unique way as $x = E(x) + \varepsilon(x)$, with $E(x) \in \mathcal{Z}$ and $|\varepsilon(x)|_\infty < 1$ (Artin decomposition; [1], theorem 12.0.3).

The uniform distribution can be defined like this ([1], definition 12.0.1):

**Definition 5.3.** A sequence $\{x_n\}_{n=1}^{\infty}$ in $\mathcal{F}_\infty$ is said to be uniformly distributed modulo $\mathcal{Z}$ if, for every $h \in \mathbb{N}$ and every $\beta \in \mathcal{F}_\infty$, we have:

$$\lim_{N\to\infty} \frac{1}{N}\, A(N; h, \beta) = q^{-h} \tag{5.10}$$

where $A(N; h, \beta)$ = number of terms $x_n$ of the series such that $n \leq N$ and $|x_n - \beta|_\infty < q^{-h}$.

Here the version of Weyl's metric theorem is the following ([1], theorem 12.0.4):

**Theorem 5.4.** *Given any $\alpha \in \mathcal{F}_\infty\ \mathcal{Z}_\infty$, the sequence $\{\lambda\,\alpha^n\}_{n=1}^{\infty}$ is uniformly distributed modulo $\mathcal{Z}$ for almost all $\lambda \in \mathcal{F}_\infty$ (in the sense of a Haar measure).*

## 6. Conclusions and future research

Most of the results presented here have been known for several decades, although I found theorems 3.5 and 4.2 independently.

Theorem 3.5 can be considered as a direct consequence of very general results in Pisot's thesis ([8]). On the other hand, the proof presented here gives more details for the specific case of a sequence of the form $\lambda\,\alpha^n$. Expression (3.8) gives a recursive rule that could be used to design an algorithm that approximates the desired value of $\alpha$.

In Pisot's thesis there are also global versions of theorem 4.2. So far, however, I have not found in any of the papers and works I have read the kind of local version given here.

Several topics that deserve further study are the following:

(1) Find a concrete (computable?) real number $\alpha > 1$ such that $\alpha^n$ is u.d. mod 1. To this end, an answer to the following question might be helpful.

(2) In section 3 we saw a method to get a number $\alpha$ such that $\lambda\,\alpha^n$ is close modulo 1 to a given sequence $\theta_n$. If the sequence $\theta_n$ is u.d. mod 1, does that imply that $\lambda\,\alpha^n$ is u.d. mod 1?

(3) Find a concrete transcendental number $\alpha > 1$ such that $\alpha^n$ is not u.d. mod 1. The results in section 3 prove that there are transcendental numbers with this property, but no explicit example is given.

(4) Is there any transcendental number $\alpha > 1$ such that $\lim_{n\to\infty} \|\alpha^n\| = 0$? According to the results in section 5.1, the convergence $\|\alpha^n\| \to 0$ should be rather slow.

(5) Generalize theorem 4 to Dedekind fields (see [7], chap. 5, for a definition of Dedekind field).

(6) Further asymptotic results: we already know (in $\mathbb{R}$) that if $\|\alpha^n - \beta^n\| = 0$ for every $n \geq n_0$ then $\alpha = \beta$ or $\alpha, \beta \in \mathbb{Z}$, but what conclusions can be drawn from $\lim_{n\to\infty} \|\alpha^n - \beta^n\| = 0$? And from $\lim_{n\to\infty} \| \sum_{i=1}^m \lambda_i \, \alpha_i^n \| = 0$? How do the results extend to local fields?

(7) Extend results to distribution over $\mathcal{O}_v$ for some ring of $v$-adic integers $\mathcal{O}_v = \{ x \in k_v \, : \, |x|_v \leq 1 \}$, to rings of adeles, and to fields of formal power series.

Department of Mathematics, RLM 8.100, University of Texas at Austin, Austin, Texas 78712

*E-mail address*: mlerma@math.utexas.edu